\documentclass[11pt]{article}
\usepackage{amsfonts}

\usepackage{amsfonts, amsmath, amssymb}
\usepackage{amssymb,amsfonts,amsmath,latexsym, epsfig,cite, psfrag,eepic,colordvi}
\usepackage{amscd,graphics}
\usepackage{lineno}

\newcommand{\qed}{\hfill $\Box $}
\newcommand{\pf}{\noindent {\bf Proof.} }

\newtheorem{theorem}{Theorem}[section]

\begin{document}

\title{Vertex-deleted subgraphs and regular factors from regular graph\thanks{This work is supported by the Fundamental Research Funds
for the Central Universities.}}

\author{Hongliang Lu\,\textsuperscript{a}\thanks{Corresponding email:
luhongliang215@sina.com (H. Lu)}, Bing Bai\,\textsuperscript{b}, Wei
Wang\,\textsuperscript{a}
\\ {\small \textsuperscript{a}Department of Mathematics}
\\ {\small Xi'an Jiaotong University, Xi'an 710049, PR China}
\\ {\small \textsuperscript{b}Center for Combinatorics, LPMC}
\\ {\small Nankai University, Tianjin, China}
}

\date{}

\maketitle

\date{}

\maketitle

\begin{abstract}
Let $k$, $m$ and $r$ be three integers such that $2\leq k\leq  m\leq
r$.
 Let $G$ be a $2r$-regular, $2m$-edge-connected graph of odd order.
We obtain some sufficient conditions, which guarantee $G-v$ contains
a $k$-factor for all $v\in V(G)$.
\end{abstract}

\section{Introduction}

All graphs considered are multigraphs (with loops) and finite. Let
$G = (V,E)$ be a graph with vertex set $V (G)$ and edge set $E(G)$.
The number of vertices of a graph $G$ is called the \emph{order} of
$G$ and is denoted by $n$. On the other hand, the number of edges of
$G$ is called the \emph{size} of $G$ and is denoted  by $e$. We
denote the degree of vertex $v$ in $G$ by  $d_{G}(v)$.
For two subsets $S,T\subseteq V(G)$, let $e_{G}(S,T)$
denote the number of edges of $G$ joining $S$ to $T$. 

Let $c(G)$ and $c_{o}(G)$ denote the number of components  and the
number of odd components of $G$, respectively.
  Let $k$ be a positive integer. A \emph{$k$-factor} of a graph
$G$ is a spanning subgraph $H$ of $G$ such that $d_{H}(x)=k$ for
every $x\in V(G)$.

Peterson obtained the following theorem, which is chronologically
the first result on $k$-factors in regular graphs.

\begin{theorem}[Petersen, \cite{Petersen}]
Every 3-regular, 2-connected graph has a 1-factor.

\end{theorem}

For 1-factor in arbitrary graphs, we have the following
characterization by Tutte \cite{Tutte47}.

\begin{theorem}[Tutte, \cite{Tutte47}]\label{Perfect matching}
A graph $G$ has a 1-factor if and only if $c_{o}(G-S)\leq |S|$ for
all $S\subseteq V(G)$.
\end{theorem}

The following theorem is well-known Tutte's $k$-factor Theorem.

\begin{theorem}[Tutte,\cite{Tutte52}]\label{Tutte}
Let $G$ be a graph and $k$ be an integer. Then $G$ has a $k$-factor
if and only if, for all $D,S\subseteq V(G)$ with $D\cap
S=\emptyset$,
\begin{align*}
\delta_{G}(D,S)=k|D|+\sum_{x\in
S}d_{G}(x)-k|S|-e_{G}(D,S)-q_{G}(D,S;k)\geq 0,
\end{align*}
where $q_{G}(D,S;k)$ is the number of  components $C$ of $G-(D\cup
S)$ such that $e_G(V(C),S)+k|C|\equiv 1\pmod 2$.  Moreover,
$\delta_{G}(D,S)\equiv k|V(G)|\pmod 2$. (Sometimes $C$ is called
$k$-odd component.)
\end{theorem}

The following theorem examines the existence of a 1-factor in
vertex-deleted subgraphs of a regular graph.

\begin{theorem}[Little et al., \cite{Little}]\label{Grant}
Let $G$ be a $2r$-regular, 2r-edge-connected graph of odd order and
$u$ be any vertex of $G$. Then the graph $G-u$ has a 1-factor.
\end{theorem}

Katerinis  presented the following result, which generalizes Theorem
\ref{Grant}.
\begin{theorem}[Katerinis,\cite{Katerinis}]\label{Katerinis}
Let $G$ be a $2r$-regular, 2r-edge-connected graph of odd order and
$m$ be an integer such that $1\leq m\leq r$. Then for every $u\in
V(G)$, the graph $G-u$ has an $m$-factor.
\end{theorem}


\section{Main result}

The purpose of this paper is to present the following result which
generalizes Theorem \ref{Katerinis}.

\begin{theorem}\label{main}
Let  $m$ and $r$ be two integers such that $2\leq m< r$.
 Let $G$ be  $2r$-regular, $2m$-edge-connected graph with odd order.
If one of the following conditions holds, then $G-v$ has a
$k$-factor for all $v\in V(G)$.
\begin{itemize}
\item[$(i)$] $k$ is even and $2\leq k\leq m$;

\item[$(ii)$] $k$ is odd, $3\leq k \leq m$ and $2m>r$.

\end{itemize}
\end{theorem}


\pf Suppose that the result doesn't hold. Then there exists $u\in
V(G)$ such that $G-u$ contains no $k$-factors. Let $H=G-u$. By
Theorem \ref{Tutte}, there exist two disjoint subsets   $D$ and $S$
of $V(G)-u$ such that
\begin{align}\label{eq:1}
q_{H}(D,S;k)+\sum_{x\in S}(k-d_{H-D}(x))\geq k|D|+2.
\end{align}
Define $S'=S\cup \{u\}$ and $W=(G-D)-S'$.

\medskip
\textbf{ Claim 1.~} $c(W)\geq 2$.
\medskip

Otherwise, suppose $c(W)\leq 1$. We consider two cases.

\medskip
{\it Case 1.} $c(W)=0$.

Since $c(W)\geq q_{H}(D,S;k)$, (\ref{eq:1}) implies
\begin{align}\label{eq:2}
\sum_{x\in S}(k-d_{H-D}(x))\geq k|D|+2.
\end{align}
So $k|S|\geq k|D|+2$ and hence $|S|>|D|$. But $V(H)=D\cup S$ and
$|V(H)|$ is even, therefore
\begin{align}\label{eq:3}
|S|\geq |D|+2.
\end{align}
Now since $G$ is $2r$-regular by Tutte's theorem we have
$$\sum_{x\in S'}(2r-d_{G-D}(x))\leq 2r|D|,$$ which implies
$$2r|S'|-\sum_{x\in S'}d_{G-D}(x)\leq 2r|D|.$$ Therefore,
$$2r(|S|+1)-\sum_{x\in S'}d_{G-D}(x)\leq 2r|D|, $$ and hence,
\begin{align}\label{eq:4}
(2r-k)|S|+k|S|+2r-\sum_{x\in S'}d_{G-D}(x)\leq k|D|+(2r-k)|D|.
\end{align}
However
\begin{align*}
\sum_{x\in S'}d_{G-D}(x)&=\sum_{x\in S}d_{G-D}(x)+d_{G-D}(u)\\
&=\sum_{x\in S}d_{H-D}(x)+e_G(u,S)+d_{G-D}(u).
\end{align*}
Therefore (\ref{eq:4}) becomes
\begin{align}\label{eq:5}
k|D|+(2r-k)|D|\geq (2r-k)|S|+k|S|+2r-\sum_{x\in
S}d_{H-D}(x)-e_G(u,S)-d_{G-D}(u)
\end{align}
Now using (\ref{eq:2}) and since $e_{G}(u,S)\leq d_{G-D}(u)\leq 2r$,
(\ref{eq:5}) implies
\begin{align}\label{eq:6}
(2r-k)(|S|-|D|)\leq 2r-2.
\end{align}
Moreover, by (\ref{eq:3}), $|S|\geq |D|+2$, we can conclude from
(\ref{eq:6}) that $k\geq r+1$. That is a contradiction, thus Case 1
can't occur.

\medskip
{\it Case 2.} $c(W)=1$.

Then we have $q_{H}(D,S;k)\leq 1$ and this implies
\begin{align}\label{eq:7}
\sum_{x\in S}(k-d_{H-D}(x))\geq k|D|+1.
\end{align}
So $k|S|\geq k|D|+1$ and hence $|S|>|D|$. 

Since $G$ is a $2m$-edge-connected, $2r$-regular graph,
so we have
\begin{align*}
2r|D|&\geq e_{G}(D,V(G-D))=e_{G}(D,V(W))+e_{G}(D,S')\\
&=e_{G}(D\cup S',V(W))-e_{G}(S',V(W))+e_{G}(D,S')\\
&=e_{G}(D\cup S',V(W))-(\sum_{x\in
S'}d_{G-D}(x)-2e_{G}(S',S'))+(2r|S'|-\sum_{x\in S'}d_{G-D}(x))\\
&\geq e_{G}(D\cup S',V(W))-2\sum_{x\in
S'}d_{G-D}(x)+2e_{G}(S',S')+2r|S'|\\
&\geq 2m-2(\sum_{x\in
S}d_{H-D}(x)+e_{G}(u,S)+d_{G-D}(u))+2e_{G}(S',S')+2r|S'|\\
&\geq 2m-2\sum_{x\in S}d_{H-D}(x)-2d_{G-D}(u)+2e_{G}(S,S)+2r|S'|\\
&\geq 2m-2\sum_{x\in S}d_{H-D}(x)+2r|S|-2r.
\end{align*}
Now using (\ref{eq:7}) implies
\begin{align*}
2r|D|&\geq 2m-2\sum_{x\in S}d_{H-D}(x)-2r+2r|S|\\
&\geq 2m-2(k|S|-k|D|-1)-2r+2r|S|.\\
\end{align*}
Thus $$(2r-2k)(|D|-|S|)\geq 2m- 2r+2\geq 2k-2r+2$$ from which it
follows $|D|\geq |S|$, a contradiction. So Case 2 can't also occur.

So we have $c(W)\geq 2$. We denote the components of $W$ by
$C_{1},\ldots,C_{c(W)}$. Suppose that $e_{G}(C_{1},D\cup
S')\leq\cdots \leq e_{G}(C_{c(W)},D\cup S')$.

Firstly, we consider (i). Then we have
\begin{align*}
2r|D|&\geq e_{G}(D\cup S',V(W))-2\sum_{x\in
S'}d_{G-D}(x)+2e_{G}(S',S')+2r|S'|\\
&\geq (2m-2)c(W)-4r+(2r-2k)|S|+2r+2k|D|+4\\
&\geq 4m-2r+(2r-2k)|S|+2k|D|.
\end{align*}
Thus we have $(2r-2k)(|D|-|S|+1)\geq 2m$ from which it follows
$|D|\geq|S|$. For every odd component $C$ of $W$, the integer
$k|V(C)|+e_{H}(V(C),S)$ is odd and since $k$ is an even integer, so
$e_{H}(V(C),S)$ must be odd. Thus $e_{H}(V(C),S)\geq 1$ and so
$\sum_{x\in S}d_{H-D}(x)\geq q_{H}(D,S;k)$. Hence (\ref{eq:1})
implies
$k|S|\geq k|D|+2$ from which it follows that $|S|\geq |D|+1$, a contradiction. 

Secondly, we show that (ii), that is, $k\geq 3$ is odd and $2m>r$.
By Theorem \ref{Katerinis}, we can assume that  $m<r$. So we have
\begin{align*}
2r|D|&\geq e_{G}(D\cup S',V(W))-2\sum_{x\in
S'}d_{G-D}(x)+2e_{G}(S',S')+2r|S'|\\
&\geq (2m-2)c(W)-4r+(2r-2k)|S|+2r+2k|D|+4\\
&\geq 4m-2r+(2r-2k)|S|+2k|D|\\& \geq (2r-2k)|S|+2k|D|+2.
\end{align*}
Thus we have $(2r-2k)(|D|-|S|)\geq 4m-2r\geq 2$ and hence $|D|>|S|$.
Let $q=q_{H}(D,S;k)$. Note that
\begin{align*}
2r|D|\geq 2mq+2r|S|-2r-2\sum_{x\in S}d_{H-D}(x).
\end{align*}
So we obtain
\begin{align}\label{eq:10}
|D|-|S|\geq \frac{m}{r}q-1-\frac{1}{r}\sum_{x\in S}d_{H-D}(x).
\end{align}
By (\ref{eq:1}), we have
\begin{align}\label{eq:11}
|D|-|S|\leq \frac{1}{k}(q-\sum_{x\in S}d_{H-D}(x)-2),
\end{align}
and $q\geq k+2$ since $|D|>|S|$. By (\ref{eq:10}) and (\ref{eq:11}),
we have
\begin{align*}
0\leq(\frac{1}{k}-\frac{1}{r})\sum_{x\in S}d_{H-D}(x)&\leq
\frac{q}{k}-\frac{2}{k}-\frac{mq}{r}+1\\
&< \frac{q}{k}-\frac{2}{k}-\frac{q}{2}+1 \\
&\leq q(\frac{1}{k}-\frac{1}{2})-\frac{2}{k}+1\\
&\leq (k+2)(\frac{1}{k}-\frac{1}{2})-\frac{2}{k}+1\\
&=1-k/2<0,
\end{align*}
a contradiction. We complete the proof. \qed

The bounds are sharp. Firstly, we show that the upper bound is
sharp. Let $G_1$ be the complete graph $K_{2r+1}$ from which a
matching of size $m$ is deleted. Let $G_{2}$ be the bipartite graph
with bipartition $(U,W)$ obtained by deleting a matching of size $m$
from $K_{2r,2r}$. Let $G$ be the $2r$-regular graph obtained by
matching the $2m$ vertices of degree $2r-1$ in $G_{1}$ to $2m$
vertices of degree $2r-1$ in $G_2$. Clearly, $G$ is
$2m$-edge-connected. Let $m^*\geq m+1$. Now we show that $G-u$
contains no $m^*$-factors for all $u\in U\cup W$. Without loss
generality, suppose that $u\in U$. Let $D=U-u$, $S=W$ and $G'=G-u$.
Note that $q_{G'}(D,S;m^*)=1$ if $m^*\neq m$ (mod 2) and
$q_{G'}(D,S;m^*)=0$ if $m^*\equiv m$ (mod 2). Then we have
\begin{align*}
m^*|D|-m^*|S|+\sum_{x\in S}d_{G'-D}(x)-q_{G'}(D,S;m^*)\leq-2<0.
\end{align*}
So by Theorem \ref{Tutte}, $G-u$ contains no $m^*$-factors.

Next we show that the lower bound is sharp. Let $\Gamma$ be the
complete graph $K_{2r+1}$ from which a matching of size $r-1$ is
deleted. Take $2r-1$ disjoint copies of $\Gamma$. Let $M_{r-1}$ be a
matching of size $r-1$. Let $H$ be the $2r$-regular graph obtained
by matching the $2r-2$ vertices of degree $2r-1$ in each of $2r-1$
disjoint copies of $\Gamma$ to the vertex set $S$ of $M_{r-1}$.
Since $c_{o}(H-V(M_{r-1}))=2r-1>|V(M_{r-1})-v|=2r-3$ for all $v\in
V(M_{r-1})$, so by Theorem \ref{Perfect matching}, $G-v$ contains no
1-factor for all $v\in V(M_{r-1})$. So the lower bound is sharp.

Finally, we show that the condition $2m>r$ is sharp. Otherwise,
suppose that $2m\leq r$. Let $R_{1}$ denote the complete bipartite
graph $K_{2r,2r-1}$ with bipartition $(U,W)$, where $|U|=2r$
 and $|W|=2r-1$. Let $R_{2}$ be the complete graph
$K_{2r+1}$ from which a matching of size $m$ is deleted. Take two
copies of $R_{2}$. Match $4m$  vertices of degree $2r-1$ of two
copies of $R_2$ to $4m$  vertices of degree $2r-1$ of $K_{2r,2r-1}$,
and then add a matching of size $r-2m$ to the rest vertices of
degree $2r-1$ of $K_{2r,2r-1}$. Then we obtain a $2r$-regular,
$2m$-edge-connected graph $R$. Let $u\in U$ and $R'=R-u$. Let
$D=U-u$ and $S=W$. Since $k$ is odd, so  $q_{R'}(D,S;k)=2$ and hence
we have
\begin{align*}
k|D|-k|S|+\sum_{x\in S}d_{R'-D}(x)-q_{G'}(D,S;k)=-2<0.
\end{align*}
 So by Theorem \ref{Tutte}, $G-u$ contains no $k$-factors.

\end{document}